\documentclass[11pt]{article}

\usepackage{authblk}
\usepackage[margin=1in]{geometry}
\usepackage[T1]{fontenc}
\usepackage{lmodern}
\usepackage{amsmath,amssymb,amsthm,mathtools}
\usepackage{microtype}
\usepackage[dvipsnames]{xcolor}
\usepackage{graphicx}
\usepackage{tikz}
\usetikzlibrary{calc,decorations.pathreplacing}
\usepackage[hidelinks]{hyperref}

\definecolor{RadoRed}{RGB}{177,43,48}
\definecolor{RadoBlue}{RGB}{40,83,155}

\theoremstyle{plain}
\newtheorem{theorem}{Theorem}[section]
\newtheorem{lemma}[theorem]{Lemma}
\newtheorem{corollary}[theorem]{Corollary}
\theoremstyle{definition}
\newtheorem{definition}[theorem]{Definition}

\theoremstyle{definition}
\newtheorem{example}[theorem]{Example}

\usepackage[nameinlink,noabbrev]{cleveref}
\crefname{theorem}{Theorem}{Theorems}
\crefname{lemma}{Lemma}{Lemmas}
\crefname{corollary}{Corollary}{Corollaries}
\crefname{figure}{Figure}{Figures}

\title{Exact Rado Numbers for Polynomial Difference Equations}
\author[1]{Ishan Kumthekar}
\author[2]{William Gasarch}

\affil[1]{Department of Mathematics, University of Florida,
Gainesville, FL 32611, USA\\
\texttt{ishankumthekar@ufl.edu}}

\affil[2]{Department of Computer Science, University of Maryland,
College Park, MD 20742, USA\\
\texttt{gasarch@umd.edu}}

\date{August 29, 2026}

\begin{document}
\maketitle
\begin{abstract}
We consider exact two-color Rado numbers for nonlinear equations of the form
$x-y=aP(z)$, where $P(z)$ is an integer-valued polynomial
satisfying $P(0)=0$ and $P(1)=1$. Doss, Saracino, and Vestal proved that
$\mathcal{R}_2(x-y=z^d)=2^{d+1}+1$ for $d\ge2$. Since the Polynomial Van der
Waerden number satisfies $W(z^d,2)=2^d+1$, their formula can be rewritten as
\[
\mathcal{R}_2(x-y=P(z))
=
P(P(1)+1)+W(P(z),2)
\]
when $P(z)=z^d$. This suggests a structural connection between Polynomial Van
der Waerden and Rado numbers. Motivated by this identity, we ask for which
integer-valued polynomials $P(z)$ and positive integers $a$ the analogous formula
\[
\mathcal{R}_2(x-y=aP(z))
=
aP(a+1)+W(aP(z),2)
\]
continues to hold. We give general conditions answering this question and show
that they extend the earlier result of Doss, Saracino, and Vestal while producing
several new infinite families of exact nonlinear Rado numbers. These include
exact formulas for equations such as    $x-y=c(z^m+z^n)$ and $x-y=a\sum_{r=1}^{z} r^d$, along with families arising from
polygonal numbers, simplex numbers, and rising factorials.
\end{abstract}

\section{Introduction}

Schur's theorem states that every finite coloring of the positive integers contains monochromatic positive integers $x,y,z$ satisfying
\[
x+y=z.
\]
The Schur number $S(r)$ is the largest $N$ for which $[N]=\{1,\ldots,N\}$ has an $r$-coloring with no monochromatic solution. The first values are
\[
S(1)=1,\qquad S(2)=4,\qquad S(3)=13.
\]
Rado extended Schur's theorem to homogeneous linear equations and systems \cite{Rado1933}. The finite form leads to Rado numbers.

\begin{definition}
Let $E$ be an equation and let $c\geq 1$. The \emph{$c$-color Rado number} of $E$, denoted by
\[
\mathcal{R}_c(E),
\]
is the least positive integer $N$ such that every $c$-coloring of $[N]$ contains a monochromatic solution to $E$. Variables in the solution need not be distinct.
\end{definition}
\begin{example}
Some two-color Rado numbers for linear equations, due to
Beutelspacher and Brestovansky \cite{BeutelspacherBrestovansky}, are
\[
\begin{aligned}
\mathcal{R}_2(x_1+\cdots+x_{m-1}=x_m)
    &=m^2-m-1, \qquad (m\ge3),\\
\mathcal{R}_2(x_1+x_2+x_3=x_4)
    &=11.
\end{aligned}
\]
Other examples with nonlinear equations are given by Doss, Saracino, and Vestal
\cite{DossSaracinoVestal}. After rearranging variables,
\[
\begin{aligned}
\mathcal{R}_2(x-y=z^d)
    &=2^{d+1}+1, \qquad (d\ge2),\\
\mathcal{R}_2(ax-y=z^2)
    &=a-1, \qquad (a\ge2).
\end{aligned}
\]
\end{example}
\begin{definition}
Let $P(z)$ take positive integer values at positive integers. The
\emph{$c$-color Polynomial Van der Waerden number},
\[
W(P(z),c),
\]
is the least positive integer $N$ such that every $c$-coloring of $[N]$
contains integers $u<v$ colored the same and satisfying
\[
v-u=P(t)
\]
for some positive integer $t$.
\end{definition}

\begin{example}
For $P(z)=z^d$ with $d\geq 2$,
\[
W(z^d,2)=2^d+1.
\]
To show $W(z^d,2)>2^d$, color $[2^d]$ by parity.
For $z=1$, we have $z^d=1$, but adjacent integers have different colors. For $z\ge2$, we have $z^d\ge2^d$ while no two integers in $[2^d]$ are at distance $2^d$ or more. Hence
there is no monochromatic $d$th-power difference in $[2^d]$.

For the upper bound, consider any $2$-coloring of $[2^d+1]$. If two
adjacent integers have the same color, they differ by $1=1^d$. Otherwise
the coloring alternates, so $1$ and $2^d+1$ have the same color. They
differ by $2^d$, giving a monochromatic $d$th-power difference. Therefore, $W(z^d,2)\leq2^d+1$.

For $P(z)=z^d$, we have $P(1)=1$. Thus, the Doss--Saracino--Vestal formula
\[
\mathcal{R}_2(x-y=z^d)=2^{d+1}+1 = 2^d + 1 + 2^d
\]
can be written as
\[
\mathcal{R}_2(x-y=P(z))=P(P(1)+1)+W(P(z),2).
\]
\end{example}
This raises the natural question of when such an identity persists beyond the
power polynomials $P(z)=z^d$.

\section{Polynomial Van der Waerden numbers}

Assume from now on that $P(z)$ is integer-valued,
\begin{equation}\label{eq:normalization}
P(0)=0,\qquad P(1)=1,
\end{equation}
and
\begin{equation}\label{eq:increasing}
1=P(1)<P(2)<P(3)<\cdots.
\end{equation}
Let
\[
h=\min\{j\ge2:P(j)\text{ is even}\}.
\]
Thus $P(1),\ldots,P(h-1)$ are odd and $P(h)$ is even.

\begin{lemma}
If $P(z)$ has degree $d\ge1$, then $h$ is well-defined and
\[
2\le h\le d+1.
\]
\end{lemma}

\begin{proof}
Since $P(1)=1$ is odd, we have $h\ge2$. Suppose, for contradiction, that
$P(2),\ldots,P(d+1)$ are all odd. Together with $P(0)=0$ and $P(1)=1$,
the parities of $P(0),P(1),\ldots,P(d+1)$ are therefore
$0,1,1,\ldots,1$.

Recall that $\Delta P(n)=P(n+1)-P(n)$. Starting from this parity sequence,
the first forward differences modulo $2$ are $1,0,0,\ldots,0$.
Taking differences again gives the same pattern, with one fewer entry,
since $0-1\equiv1\pmod2$ and all remaining differences are $0$.
Repeating this process, after $d+1$ differences only the $1$ remains.
Hence
\[
\Delta^{d+1}P(0)\equiv1\pmod2.
\]
But each forward difference lowers the degree by one, so the $(d+1)$st
forward difference of a degree-$d$ polynomial is identically zero. This is
a contradiction. Therefore some $P(j)$ with $2\le j\le d+1$ is even, and
hence $h\le d+1$.
\end{proof}

\begin{lemma}\label{lem:first-even}
Under \eqref{eq:normalization} and \eqref{eq:increasing},
\[
W(P(z),2)=P(h)+1.
\]
\end{lemma}

\begin{proof}
Color $[P(h)]$ by parity. If $j<h$, then $P(j)$ is odd, so two integers differing by $P(j)$ have opposite parity and therefore cannot have the same color. If $j\ge h$, then $P(j)\ge P(h)$ and cannot occur inside $[P(h)]$. Hence $W(P(z),2)>P(h)$.

Conversely, consider any coloring of $[P(h)+1]$. If two consecutive integers have the same color, their difference is $P(1)=1$. Otherwise the colors alternate. Since $P(h)$ is even, $1$ and $P(h)+1$ then have the same color and differ by $P(h)$. Thus $W(P(z),2)\le P(h)+1$.
\end{proof}

\begin{lemma}\label{lem:scaling}
For every positive integer $a$,
\[
W(aP(z),2)=a\bigl(W(P(z),2)-1\bigr)+1
\]
Consequently,
\[
W(aP(z),2)=aP(h)+1
\]
\end{lemma}
Gasarch, Kruskal, Kruskal, and Price prove the first identity for
$P\in\mathbb Z[z]$ \cite{GasarchKruskalKruskalPrice}.
Their proof applies unchanged to integer-valued polynomials. The second
identity follows from Lemma~\ref{lem:first-even}.
\section{Main theorem}

\begin{theorem}\label{thm:main}
Let $P(z)$ be an integer-valued polynomial satisfying
\eqref{eq:normalization} and \eqref{eq:increasing}, and let
\[
h=\min\{j\ge2:P(j)\text{ is even}\}.
\]
Let $a\ge h$, and assume
\begin{equation}\label{eq:range}
P(2a+1)\ge P(a+1)+P(h).
\end{equation}
Then
\begin{equation}\label{eq:main-identity}
\boxed{
\mathcal{R}_2(x-y=aP(z))
=
aP(a+1)+W(aP(z),2)
}
\end{equation}
if and only if there is no $j$ with $1\le j\le 2a$ such that
\begin{equation}\label{eq:criterion}
0<|P(j)-P(a+1)|<P(h)-1
\quad\text{and}\quad
|P(j)-P(a+1)|\text{ is odd}.
\end{equation}
\end{theorem}

We first prove the upper bound.

\begin{lemma}\label{lem:upper}
Under the assumptions of \cref{thm:main},
\[
\mathcal{R}_2(x-y=aP(z))
\le
aP(a+1)+aP(h)+1.
\]
\end{lemma}

\begin{proof}
The following upper bound argument generalizes the construction
of Doss, Saracino, and Vestal \cite{DossSaracinoVestal}. Their argument uses
the distances $1$ and $2^d$, while here we replace these with general distances $aP(1)=a$ and $aP(h)$.
Set $N=aP(a+1)+aP(h)+1$ and suppose that $[N]$ has no monochromatic solution. We first show that $1$ and $h$ have the same color. After interchanging colors if necessary, suppose that $1$ is red and $h$ is blue, and put $H=P(h)$. Consider
\[
u_i=1+ia,\qquad 0\le i\le2H.
\]

These all lie in $[N]$, since $H<P(a+1)$. No two consecutive $u_i$ can both be red, since their difference is $a=aP(1)$. Also, if $u_i$ is blue for $0\le i\le H$, then $u_{i+H}$ must be red, since their difference is $aH=aP(h)$.

Now $u_0$ is red, so $u_1$ is blue, hence $u_{H+1}$ is red. Thus $u_H$ is blue, and hence $u_{2H}$ is red. Moreover, for each $1\le i\le H-1$, at least one of $u_i$ and $u_{i+H}$ is red. Therefore there are at least
\[
(H-1)+2=H+1
\]
red vertices among $u_0,\ldots,u_{2H}$. Since no two consecutive vertices are red, this is the maximum possible number, forcing the red vertices to be
\[
u_0,u_2,u_4,\ldots,u_{2H}.
\]
But $H=P(h)$ is even, so $u_H$ is red, contradicting that $u_H$ is blue. Hence $1$ and $h$ have the same color. Without loss of generality assume they are red.

Suppose now that some $z\le a+1$ is blue. Consider the two rows
\[
\resizebox{0.98\linewidth}{!}{$\displaystyle
\begin{array}{c|cccccc}
\text{column}
&0&1&2&\cdots&P(h)-1&P(h)\\
\hline
\text{top}
&1&1+a&1+2a&\cdots&1+(P(h)-1)a&1+aP(h)\\
\text{bottom}
&1+aP(z)&1+aP(z)+a&1+aP(z)+2a
&\cdots
&1+aP(z)+(P(h)-1)a
&1+aP(z)+aP(h)
\end{array}
$}
\]

Each row is a copy of the odd cycle whose consecutive edges have length
$aP(1)=a$ and whose closing edge has length $aP(h)$. The two entries in each column differ by $aP(z)$, so they
cannot both be blue. Hence each column contains a red entry.

There are now two possibilities.
\begin{enumerate}
    \item Two adjacent columns contain red entries in the same row. Then
those two red entries differ by $a=aP(1)$, so together with the red root
$1$ they give a red solution.
    \item Otherwise each column has exactly one red entry, and these red entries alternate between the two rows.
Since $P(h)$ is even, columns $0$ and $P(h)$ contain red entries in the
same row. These two entries differ by $aP(h)$, so together with the red
root $h$ they give a red solution.
\end{enumerate}
Both cases are impossible. Therefore every $z\le a+1$ is red. In particular, $1$ and $a+1$ are red,
but
\[
(a+1)-1=a=aP(1),
\]
so $(a+1,1,1)$
is a red solution, a contradiction.
\end{proof}
\vspace{1.5cm}

\begin{figure}[htbp]
\centering
\begin{tikzpicture}[
    x=1cm,
    y=1cm,
    every node/.style={font=\small}
]

\path[use as bounding box] (-0.8,-1.75) rectangle (13.1,3.45);


\node[circle,draw=gray!65,fill=RadoRed,text=white,
      minimum size=7.5mm,inner sep=0pt]
      (Ltop0) at (0,1.35) {\bfseries R};

\node[circle,draw=gray!65,fill=RadoRed,text=white,
      minimum size=7.5mm,inner sep=0pt]
      (Ltop1) at (2.35,1.35) {\bfseries R};

\node[circle,draw=gray!65,fill=white,
      minimum size=7.5mm,inner sep=0pt]
      (Lbot0) at (0,-0.10) {};

\node[circle,draw=gray!65,fill=white,
      minimum size=7.5mm,inner sep=0pt]
      (Lbot1) at (2.35,-0.10) {};

\draw[gray!65,thick] (Ltop0)--(Ltop1);
\draw[gray!65,thick] (Lbot0)--(Lbot1);

\draw[dashed,gray!65]
      (Ltop0)--(Lbot0);

\draw[dashed,gray!65]
      (Ltop1)--(Lbot1);

\draw[<->,gray!70]
      (0,1.88)--(2.35,1.88)
      node[midway,above=2pt] {$a$};

\node[left=9pt]
      at ($(Ltop0)!0.5!(Lbot0)$)
      {$aP(z)$};

\node at (0,-0.70) {$j$};
\node at (2.35,-0.70) {$j+1$};

\node[font=\bfseries]
      at (1.175,-1.40)
      {Adjacent columns};

\begin{scope}[shift={(5.0,0)}]

\node[circle,draw=gray!65,fill=RadoRed,text=white,
      minimum size=7.5mm,inner sep=0pt]
      (T0) at (0,1.35) {\bfseries R};

\node[circle,draw=gray!65,fill=white,
      minimum size=7.5mm,inner sep=0pt]
      (T1) at (2.35,1.35) {};

\node[circle,draw=gray!65,fill=RadoRed,text=white,
      minimum size=7.5mm,inner sep=0pt]
      (T2) at (4.70,1.35) {\bfseries R};

\node at (6.10,1.35) {$\cdots$};

\node[circle,draw=gray!65,fill=white,
      minimum size=7.5mm,inner sep=0pt]
      (Tm) at (7.50,1.35) {};

\node[circle,draw=gray!65,fill=RadoRed,text=white,
      minimum size=7.5mm,inner sep=0pt]
      (Th) at (9.85,1.35) {\bfseries R};

\node[circle,draw=gray!65,fill=white,
      minimum size=7.5mm,inner sep=0pt]
      (B0) at (0,-0.10) {};

\node[circle,draw=gray!65,fill=RadoRed,text=white,
      minimum size=7.5mm,inner sep=0pt]
      (B1) at (2.35,-0.10) {\bfseries R};

\node[circle,draw=gray!65,fill=white,
      minimum size=7.5mm,inner sep=0pt]
      (B2) at (4.70,-0.10) {};

\node at (6.10,-0.10) {$\cdots$};

\node[circle,draw=gray!65,fill=RadoRed,text=white,
      minimum size=7.5mm,inner sep=0pt]
      (Bm) at (7.50,-0.10) {\bfseries R};

\node[circle,draw=gray!65,fill=white,
      minimum size=7.5mm,inner sep=0pt]
      (Bh) at (9.85,-0.10) {};

\draw[gray!65,thick] (T0)--(T1)--(T2);
\draw[gray!65,thick] (Tm)--(Th);

\draw[gray!65,thick] (B0)--(B1)--(B2);
\draw[gray!65,thick] (Bm)--(Bh);

\draw[dashed,gray!65] (T0)--(B0);
\draw[dashed,gray!65] (T1)--(B1);
\draw[dashed,gray!65] (T2)--(B2);
\draw[dashed,gray!65] (Tm)--(Bm);
\draw[dashed,gray!65] (Th)--(Bh);

\draw[RadoRed,thick]
      (T0)
      to[out=68,in=112]
      node[midway,above=9pt,text=black]
      {closing edge of length $aP(h)$}
      (Th);

\node[left=9pt]
      at ($(T0)!0.5!(B0)$)
      {$aP(z)$};

\node at (0,-0.70) {$0$};
\node at (2.35,-0.70) {$1$};
\node at (4.70,-0.70) {$2$};

\node at (7.50,-0.70) {$P(h)-1$};
\node at (9.85,-0.70) {$P(h)$};

\node[font=\bfseries]
      at (4.925,-1.40)
      {Alternating columns};

\end{scope}

\end{tikzpicture}%
\caption{The two possible arrangements of the selected red entries.}
\label{fig:upper-two-cases}
\end{figure}
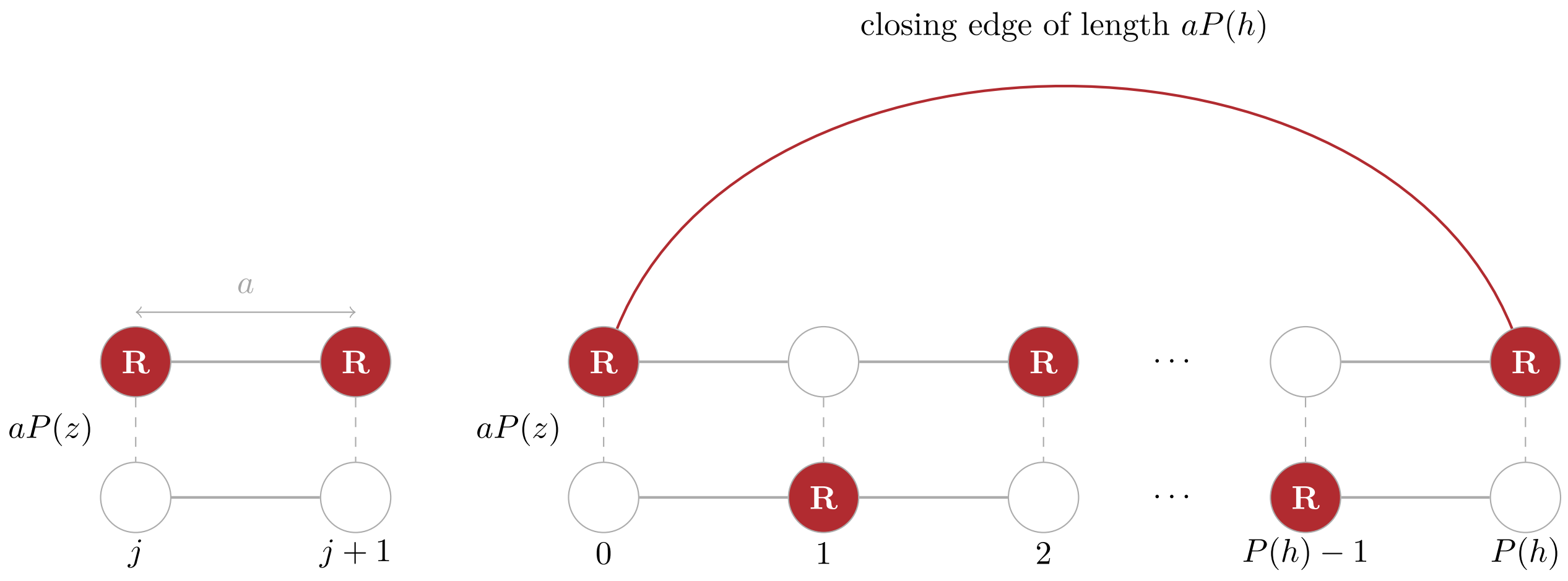

\begin{lemma}\label{lem:forced}
Under the assumptions of \cref{thm:main}, suppose
\[
[a(P(a+1)+P(h))]
\]
has no monochromatic solution. After interchanging colors if necessary,
\[
1,\ldots,a\text{ are red},
\qquad
a+1,\ldots,2a\text{ are blue}.
\]
Moreover, for $0\le i\le P(h)-1$,
\[
1+ia\text{ is red}\iff i\text{ is even},
\]
\[
1+aP(a+1)+ia\text{ is red}\iff i\text{ is odd},
\]
and $1+aP(h)$ is blue.
\end{lemma}

\begin{proof}
Using the proof of Lemma~\ref{lem:upper}, the roots $1$ and $h$ have the same color since
$1+2aP(h)\le a(P(a+1)+P(h))$.
Assume they are red. If some $z\le a$ were blue, the argument
from Lemma~\ref{lem:upper} would give a red solution. The construction
still fits inside the present interval because
$1+aP(h)+aP(z)\le 1+aP(h)+aP(a)\le a(P(a+1)+P(h))$.
Therefore $1,\ldots,a$ are red.

Now let $a+1\le z\le2a$. Then $z-a\in[1,a]$, so $z-a$ and $1$ are
red. If $z$ were also red, then $
z-(z-a)=a=aP(1),
$
and $(z,z-a,1)$ would be a red solution. Therefore
$a+1,\ldots,2a$ are blue.

Consider the two rows
\[
\begin{array}{c|ccccc}
\text{column}
&0&1&\cdots&P(h)-1&P(h)\\
\hline
\text{upper row}
&1+aP(a+1)
&1+aP(a+1)+a
&\cdots
&1+aP(a+1)+a(P(h)-1)
&\\
\text{lower row}
&1
&1+a
&\cdots
&1+a(P(h)-1)
&1+aP(h)
\end{array}
\]
The lower row has $P(h)+1$ vertices. Consecutive vertices differ by
$a=aP(1)$, while its first and last vertices differ by $aP(h)$.
These edges form an odd cycle. Since the roots $1$ and $h$ are red,
no edge of this cycle can have two red endpoints. Blue-blue endpoints are allowed, since the roots $1$ and $h$
corresponding to these edge lengths are red. Since $P(h)$ is even,
at most $\lfloor \frac{P(h)+1}{2}\rfloor= \frac{P(h)}{2}$ vertices in the lower row are red.
The upper row consists of $P(h)$ vertices whose consecutive vertices
also differ by $a=aP(1)$. Since root $1$ is red, no two consecutive
vertices can both be red. Thus at most every other vertex can be red,
so, since $P(h)$ is even, the upper row contains at most $\frac{P(h)}2$ red vertices.

Now look at the vertical columns. For every
$0\le i\le P(h)-1$, the two entries in column $i$ differ by
$aP(a+1)$. Since the root $a+1$ is blue, the two entries in a column
cannot both be blue. Thus each of the $P(h)$ columns contains at least
one red vertex.

The two rows contain at most $P(h)$ red vertices altogether, while the
$P(h)$ columns each require at least one red vertex. Therefore every
column contains exactly one red vertex, the two entries in each column
have opposite colors, and the unpaired final vertex $1+aP(h)$ is blue.

For $0\le i\le P(h)-2$, two consecutive vertices in the lower row
cannot both be red, since they differ by $a=aP(1)$ and the root $1$
is red. They also cannot both be blue, since then the two vertices
directly above them would both be red and would differ by $a=aP(1)$.
Hence the first $P(h)$ vertices of the lower row alternate in color.

Its first vertex is $1$, which is red. Therefore
\[
1+ia\text{ is red}\iff i\text{ is even}
\qquad(0\le i\le P(h)-1).
\]
The unpaired blue vertex is not part of this alternation, so the final
horizontal edge may have two blue endpoints. The upper row has the
opposite colors, proving the lemma.
\end{proof}

\begin{figure}[htbp]
\centering
\begin{tikzpicture}[x=1cm,y=1cm,every node/.style={font=\small}]

\def\yupper{2.1}
\def\ylower{0}

\node[anchor=east,font=\small\bfseries]
    at (-0.65,\yupper) {upper row};

\node[anchor=east,font=\small\bfseries]
    at (-0.65,\ylower) {lower row};


\node[
    circle,
    draw=gray!65,
    fill=RadoBlue,
    text=white,
    minimum size=7.5mm,
    inner sep=0pt
] (u0) at (0,\yupper) {\bfseries B};

\node[
    circle,
    draw=gray!65,
    fill=RadoRed,
    text=white,
    minimum size=7.5mm,
    inner sep=0pt
] (u1) at (3.0,\yupper) {\bfseries R};

\node at (6.0,\yupper) {$\cdots$};

\node[
    circle,
    draw=gray!65,
    fill=RadoRed,
    text=white,
    minimum size=7.5mm,
    inner sep=0pt
] (um) at (9.1,\yupper) {\bfseries R};


\node[
    circle,
    draw=gray!65,
    fill=RadoRed,
    text=white,
    minimum size=7.5mm,
    inner sep=0pt
] (l0) at (0,\ylower) {\bfseries R};

\node[
    circle,
    draw=gray!65,
    fill=RadoBlue,
    text=white,
    minimum size=7.5mm,
    inner sep=0pt
] (l1) at (3.0,\ylower) {\bfseries B};

\node at (6.0,\ylower) {$\cdots$};

\node[
    circle,
    draw=gray!65,
    fill=RadoBlue,
    text=white,
    minimum size=7.5mm,
    inner sep=0pt
] (lm) at (9.1,\ylower) {\bfseries B};

\node[
    circle,
    draw=gray!65,
    fill=RadoBlue,
    text=white,
    minimum size=7.5mm,
    inner sep=0pt
] (lh) at (12.2,\ylower) {\bfseries B};


\draw[gray!65,thick]
    (u0)--(u1)
    node[
        midway,
        below=5pt,
        fill=white,
        inner sep=1pt,
        font=\footnotesize,
        text=black
    ] {$a$};

\draw[gray!65,thick]
    (u1)--(5.15,\yupper);

\draw[gray!65,thick]
    (6.85,\yupper)--(um);

\draw[gray!65,thick]
    (l0)--(l1);

\draw[gray!65,thick]
    (l1)--(5.15,\ylower);

\draw[gray!65,thick]
    (6.85,\ylower)--(lm)--(lh);


\draw[dashed,gray!65] (u0)--(l0);
\draw[dashed,gray!65] (u1)--(l1);
\draw[dashed,gray!65] (um)--(lm);


\node[above=10pt,font=\footnotesize]
    at (u0)
    {$1+aP(a+1)$};

\node[above=10pt,font=\footnotesize]
    at (u1)
    {$1+aP(a+1)+a$};

\node[above=10pt,font=\footnotesize]
    at (um)
    {$1+aP(a+1)+a(P(h)-1)$};

\node[below=8pt,font=\footnotesize]
    at (l0)
    {$1$};

\node[below=8pt,font=\footnotesize]
    at (l1)
    {$1+a$};

\node[below=8pt,font=\footnotesize]
    at (lm)
    {$1+a(P(h)-1)$};

\node[below=8pt,font=\footnotesize]
    at (lh)
    {$1+aP(h)$};


\draw[<->,gray!75]
    ($(u0)+(-0.65,-0.15)$)
    --
    ($(l0)+(-0.65,0.15)$)
    node[
        midway,
        left=3pt,
        text=black,
        font=\footnotesize
    ] {$aP(a+1)$};


\draw[gray!70,thick]
    (l0)
    to[out=-40,in=-140]
    node[
        midway,
        below=18pt,
        text=black,
        font=\footnotesize
    ] {closing edge of length $aP(h)$}
    (lh);

\end{tikzpicture}%

\caption{The lower row becomes an odd cycle by adding the edge of length $aP(h)$.}
\label{fig:forcedpattern}
\end{figure}
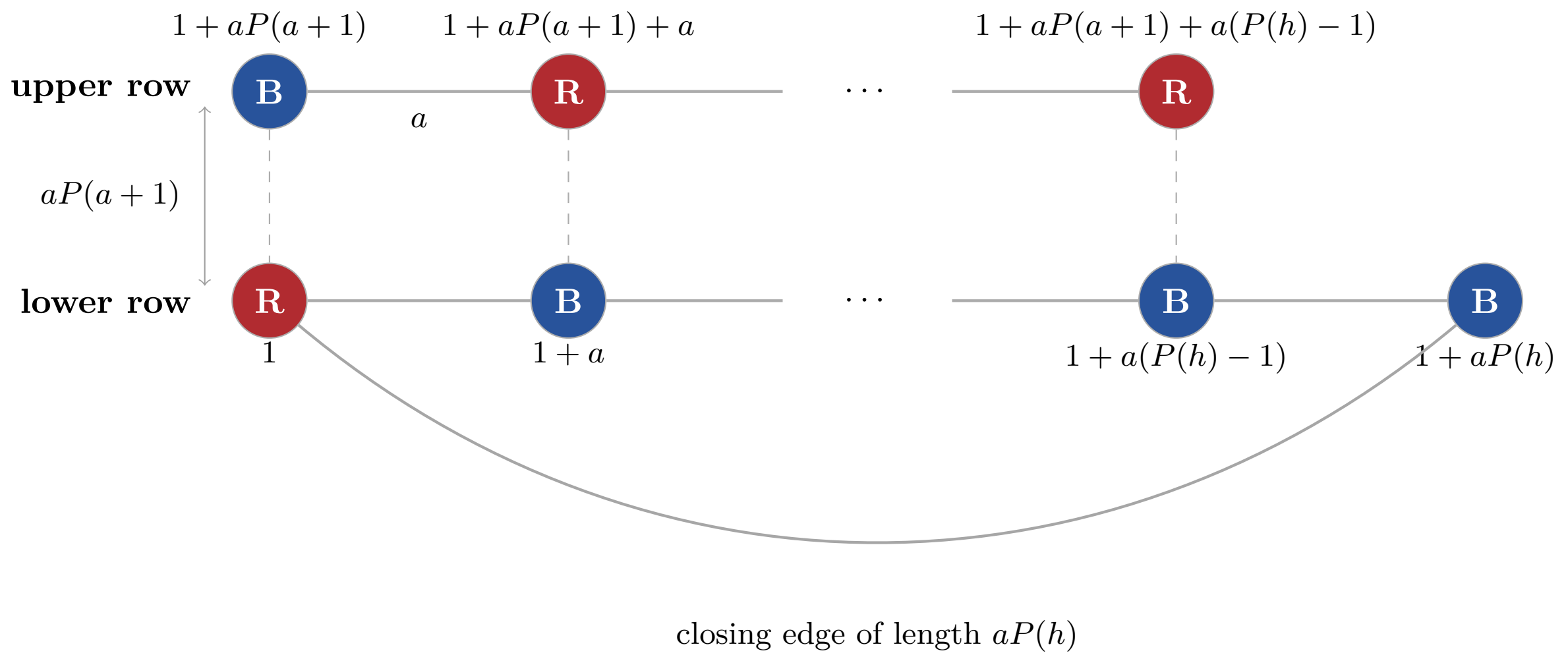

\begin{proof}[Proof of \cref{thm:main}]
Put
\[
N=a\bigl(P(a+1)+P(h)\bigr)+1.
\]
By Lemmas~\ref{lem:first-even} and~\ref{lem:scaling},
$W(aP(z),2)=aP(h)+1$. Hence the right side of
\eqref{eq:main-identity} is exactly $N$, and Lemma~\ref{lem:upper}
gives
\[
\mathcal{R}_2(x-y=aP(z))\le N.
\]
Thus the question is whether this upper bound is sharp. By the
definition of the Rado number,
\[
\mathcal{R}_2(x-y=aP(z))=N
\]
if and only if $[N-1]$ has a coloring with no monochromatic solution.
We prove that such a coloring exists if and only if no
$j$ satisfies \eqref{eq:criterion}.

First suppose that $[N-1]$ has a coloring with no monochromatic
solution. By Lemma~\ref{lem:forced}, after interchanging colors if
necessary,
$1,\ldots,a$ are red and $a+1,\ldots,2a$ are blue. Moreover,
for $0\le i\le P(h)-1$,
\[
1+ia\text{ is red}\iff i\text{ is even},
\qquad
1+aP(a+1)+ia\text{ is red}\iff i\text{ is odd}.
\]
We show that no $j$ satisfies \eqref{eq:criterion}. Suppose instead
that some $1\le j\le2a$ does. Then
\[
r=\bigl|P(j)-P(a+1)\bigr|
\]
is positive and odd, with $r<P(h)-1$. Since $P(h)$ is even, this means
$r\le P(h)-3$.

\begin{enumerate}
\item
Suppose $1\le j\le a$. Since $P$ is increasing,
$r=P(a+1)-P(j)$. The number $r+1$ is even, so
$1+a(r+1)$ is red. Also $1+aP(a+1)+a$ is red, since it corresponds
to the odd index $1$ in the upper row of Lemma~\ref{lem:forced}.
Their difference is
\[
\bigl(1+aP(a+1)+a\bigr)-\bigl(1+a(r+1)\bigr)
=a\bigl(P(a+1)-r\bigr)
=aP(j).
\]
The root $j$ is also red because $j\le a$. Hence these three integers
give a red solution, a contradiction.

\item
Suppose $a+1\le j\le2a$. Since $r>0$, we have $j>a+1$, and therefore
$r=P(j)-P(a+1)$. The number $1+a$ is blue. Since $r+1$ is even,
$1+aP(a+1)+a(r+1)$ is also blue. Their difference is
\[
\bigl(1+aP(a+1)+a(r+1)\bigr)-(1+a)
=a\bigl(P(a+1)+r\bigr)
=aP(j).
\]
The root $j$ is blue because $a+1\le j\le2a$. Hence these three
integers give a blue solution, again a contradiction.
\end{enumerate}

Therefore the absence of any $j$ satisfying \eqref{eq:criterion} is
necessary.

We now prove the converse. Assume that no $j$ satisfies
\eqref{eq:criterion}. We construct
a coloring of $[N-1]$ with no monochromatic solution. This construction
extends the alternating end coloring used in \cite{DossSaracinoVestal}.

Divide $[N-1]$ into consecutive blocks of length $a$:
\[
B_i=\{(i-1)a+1,\ldots,ia\},
\qquad
1\le i\le P(a+1)+P(h).
\]
Color the blocks as follows. For $1\le i\le P(h)$, color $B_i$ red
when $i$ is odd and blue when $i$ is even. Color
$B_{P(h)+1},\ldots,B_{P(a+1)}$ blue. Finally, for
$1\le i\le P(h)$, color $B_{P(a+1)+i}$ blue when $i$ is odd and red
when $i$ is even. This coloring is shown in
\cref{fig:block-coloring-main}.

\begin{figure}[htbp]
\centering
{
\definecolor{blockRed}{HTML}{793F3A}
\definecolor{blockBlue}{HTML}{43596D}
\begin{tikzpicture}[x=.75cm,y=1cm,
    every node/.style={inner sep=0pt},
    letter/.style={font=\small,text=white},
    count/.style={font=\footnotesize,anchor=north}]
  \foreach \start/\stop/\ink/\symbol in {
    0/1/blockRed/R,1/2/blockBlue/B,
    3/4/blockRed/R,4/5/blockBlue/B,
    5.25/9.75/blockBlue/B,
    10/11/blockBlue/B,11/12/blockRed/R,
    13/14/blockBlue/B,14/15/blockRed/R}{
      \fill[\ink] (\start,0) rectangle (\stop,.55);
      \node[letter] at ({(\start+\stop)/2},.275) {\symbol};
  }
  \foreach \x in {1,4,11,14}
    \draw[white,line width=.4pt] (\x,0)--(\x,.55);
  \foreach \x in {2.5,12.5}
    \node at (\x,.275) {$\cdots$};
  \node[count] at (2.5,-.15) {$P(h)$ blocks};
  \node[count] at (7.5,-.15) {$P(a+1)-P(h)$ blocks};
  \node[count] at (12.5,-.15) {$P(h)$ blocks};
  \node[font=\footnotesize] at (7.5,-.83) {Each block has length $a$.};
\end{tikzpicture}
\par}
\caption{A coloring of $[N-1]$ with no monochromatic solution}
\label{fig:block-coloring-main}
\end{figure}
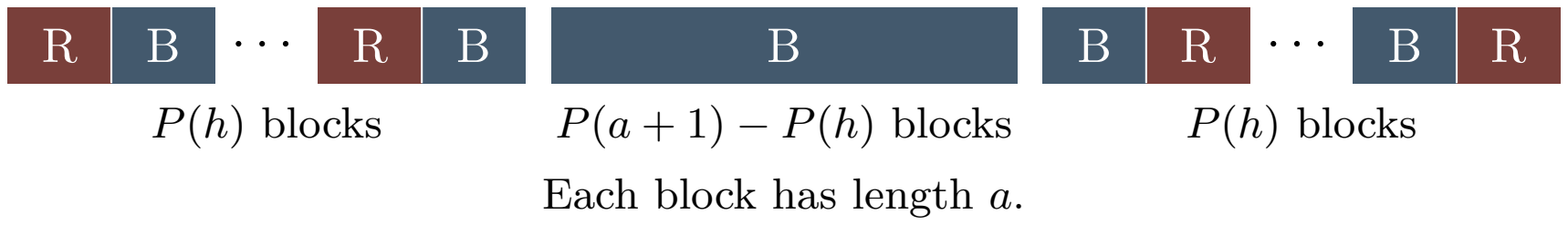
We first reduce the problem to the block indices. Suppose $x\in B_i$
and $y\in B_k$. Write
$x=(i-1)a+s$ and $y=(k-1)a+t$, where $1\le s,t\le a$. If
$x-y=aP(z)$, then
\[
a(i-k)+(s-t)=aP(z).
\]
Thus $s-t$ is divisible by $a$. Since $|s-t|<a$, we must have $s=t$,
and consequently $i-k=P(z)$.
Therefore a difference $aP(z)$ occurs exactly between equal positions
in two blocks whose indices differ by $P(z)$. To rule out a
monochromatic solution with root $z$, it is enough to show that no two
blocks having the color of $z$ have index difference $P(z)$.

We consider three possible ranges for $z$.

\begin{enumerate}
\item
Let $1\le z\le a$. Then $z\in B_1$, so the root $z$ is red. We must
show that no two red blocks have index difference $P(z)$.

First consider two red blocks lying in the same alternating end. Their
indices have the same parity, so their separation is even and strictly
less than $P(h)$. This cannot equal $P(z)$. If $z<h$, then $P(z)$ is
odd by the definition of $h$. If $z\ge h$, then
$P(z)\ge P(h)$ because $P$ is increasing.

It remains to consider a red block from the first alternating end and
a red block from the final alternating end. A red block in the first
end has index $i$, where $i$ is odd, while a red block in the final
end has index $P(a+1)+k$, where $k$ is even.

Since $z\le a$, we have $P(z)<P(a+1)$. Thus if these two blocks were
separated by $P(z)$, then
\[
P(a+1)+k-i=P(z),
\]
so
\[
P(a+1)-P(z)=i-k.
\]
The number $i-k$ is positive and odd. Moreover, since
$i\le P(h)-1$ and $k\ge2$,
\[
1\le i-k\le P(h)-3.
\]
This is precisely one of the values excluded by
\eqref{eq:criterion}. Therefore no two red blocks are separated by
$P(z)$, and no red solution can occur.

\item
Let $a+1\le z\le2a$. Then $z\in B_2$, so the root $z$ is blue. Since
$P$ is increasing, write
\[
P(z)=P(a+1)+r,
\qquad r\ge0.
\]

If $r\ge P(h)$, then $P(z)\ge P(a+1)+P(h)$, which is at least the
total number of blocks. Such a block separation cannot occur.

Suppose instead that $0\le r<P(h)$, and suppose two blue blocks are
separated by $P(z)=P(a+1)+r$. If the first block is $B_i$, then the
second must be
\[
B_{P(a+1)+r+i}.
\]
Since there are only $P(a+1)+P(h)$ blocks altogether, we must have
$i\le P(h)-r$. In particular, the first block lies in the initial
alternating end, and the second lies in the final alternating end.

In the initial end, $B_i$ is blue exactly when $i$ is even. In the
final end, $B_{P(a+1)+k}$ is blue exactly when $k$ is odd. Here
$k=r+i$. Thus, if both blocks are blue, then $i$ is even and $r+i$
is odd. Hence $r$ must be odd.

Since $i$ is a positive even integer, $i\ge2$. Together with
$i\le P(h)-r$, this gives $r\le P(h)-2$. As $r$ is odd and $P(h)$ is
even, in fact
\[
1\le r\le P(h)-3.
\]
But then
$|P(z)-P(a+1)|=r$ is a positive odd integer smaller than $P(h)-1$,
contrary to \eqref{eq:criterion}. Therefore no blue solution can occur.

\item
Let $z\ge2a+1$. By \eqref{eq:range},
\[
P(z)\ge P(2a+1)\ge P(a+1)+P(h).
\]
Hence
\[
aP(z)\ge a\bigl(P(a+1)+P(h)\bigr)=N-1.
\]
But the difference between two integers in $[N-1]$ is at most $N-2$.
Thus no pair in $[N-1]$ can have difference $aP(z)$.
\end{enumerate}
The coloring of \cref{fig:block-coloring-main} therefore contains no
monochromatic solution. Hence
$\mathcal{R}_2(x-y=aP(z))>N-1$, and so
$\mathcal{R}_2(x-y=aP(z))\ge N$. Together with the upper bound
$\mathcal{R}_2(x-y=aP(z))\le N$, we obtain
\[
\mathcal{R}_2(x-y=aP(z))
=
N
=
aP(a+1)+aP(h)+1
=
aP(a+1)+W(aP(z),2).
\]
This proves \cref{thm:main}.
\end{proof}

\section{Simplifying the hypotheses}

The exact criterion in \cref{thm:main} can be tedious to verify
separately for each polynomial family. Beyond $a\ge h$, the two
conditions to check are
\begin{enumerate}
    \item The range condition
\[
P(2a+1)\ge P(a+1)+P(h)
\]
\item The absence of a $1\le j\le2a$ for which
\[
0<|P(j)-P(a+1)|<P(h)-1
\]
and $|P(j)-P(a+1)|$ is odd.
\end{enumerate}
The next two corollaries replace these conditions by simpler hypotheses
on the consecutive gaps of $P(z)$.

\begin{corollary}\label{cor:convex}
Let $P(z)$ be an integer-valued polynomial satisfying
$P(0)=0$, $P(1)=1$, and
\[
1=P(1)<P(2)<P(3)<\cdots.
\]
Suppose that its consecutive gaps strictly increase:
\[
P(n+2)-P(n+1)>P(n+1)-P(n)
\qquad(n\ge1).
\]
Then for every $a\ge h$,
\[
P(2a+1)\ge P(a+1)+P(h).
\]
Thus the range condition in \cref{thm:main} is satisfied.
\end{corollary}

\begin{proof}
Write $\Delta_n=P(n+1)-P(n)$. Since the $\Delta_n$ are strictly
increasing integers,
$\Delta_{n+1}\ge\Delta_n+1$. Applying this $a$ times gives
\[
\Delta_{a+j}\ge\Delta_j+a.
\]
Since $a\ge h$,
\[
\begin{aligned}
P(2a+1)-P(a+1)
&\ge P(a+h)-P(a+1)\\
&=\sum_{j=1}^{h-1}\Delta_{a+j}\\
&\ge\sum_{j=1}^{h-1}(\Delta_j+a)\\
&=P(h)-1+a(h-1)\\
&\ge P(h).
\end{aligned}
\]
Hence $P(2a+1)\ge P(a+1)+P(h)$.
\end{proof}

After \cref{cor:convex}, only the forbidden odd differences in
\eqref{eq:criterion} remain to be checked. The next corollary replaces
that condition by a single inequality.

\begin{corollary}\label{cor:gap}
Let $P(z)$ be an integer-valued polynomial satisfying
$P(0)=0$, $P(1)=1$, and
\[
1=P(1)<P(2)<P(3)<\cdots,
\]
and suppose that its consecutive gaps strictly increase. If $a\ge h$
and
\[
P(a+1)-P(a)\ge P(h)-2,
\]
then
\[
\mathcal{R}_2(x-y=aP(z))
=
aP(a+1)+W(aP(z),2).
\]
\end{corollary}

\begin{proof}
By \cref{cor:convex}, the range condition already holds. It remains to verify that no $j$ satisfies \eqref{eq:criterion}.

If $j\le a$, then
\[
P(a+1)-P(j)
\ge P(a+1)-P(a)
\ge P(h)-2.
\]
Since $P(h)$ is even, every positive odd integer smaller than $P(h)-1$
is at most $P(h)-3$. Thus no $j\le a$ gives a forbidden difference.

If $j\ge a+2$, then the strictly increasing gaps give
\[
P(a+2)-P(a+1)
>
P(a+1)-P(a)
\ge P(h)-2.
\]
The gaps are integers, so
$P(a+2)-P(a+1)\ge P(h)-1$, and every later value $P(j)$ is still
farther from $P(a+1)$. Thus no $j\ge a+2$ gives a forbidden difference.
Finally, $j=a+1$ gives difference zero.

Hence no $j$ satisfies \eqref{eq:criterion}, and the result follows
from \cref{thm:main}.
\end{proof}

Thus, in the applications below, it is enough to check that $P(z)$ is
integer-valued, has $P(0)=0$ and $P(1)=1$, and is increasing.
We then determine $h$, check that the consecutive gaps strictly
increase, and verify only
\[
P(a+1)-P(a)\ge P(h)-2 \qquad a\ge h.
\]
Corollary~\ref{cor:gap} then gives the exact Rado number. These
conditions are sufficient rather than necessary. If the final gap
inequality fails, \cref{thm:main} may still apply directly, as we will
see for the simplex numbers below.

\section{Applications}

We now apply \cref{cor:gap} to get several infinite families of exact Rado numbers.

\subsection{Powers}

Let $P(z)=z^d$, where $d\ge2$. We have $h=2$ and $P(h)=2^d$.
The consecutive gaps
\[
(n+1)^d-n^d
\]
strictly increase. For every $a\ge2$, the gap at $a$ is larger than the
first gap $2^d-1$, and therefore certainly at least $P(h)-2$.
Corollary~\ref{cor:gap} gives
\[
\boxed{
\mathcal{R}_2(x-y=az^d)
=
a(a+1)^d+a2^d+1
\qquad(a\ge2).
}
\]
Notice that the case $a=1$ is the theorem of Doss, Saracino, and Vestal,
so the equation above holds for $a\ge1$.

\subsection{Sums of powers}

Let $p$ be a prime and
let
\[
d_1,\ldots,d_p\ge2,
\qquad
d_1\equiv d_2\equiv\cdots\equiv d_p\pmod{p-1}.
\]
Consider
\[
P(z)=\frac{1}{p}\sum_{i=1}^{p}z^{d_i}.
\]
This polynomial is integer-valued. Indeed, if $p\mid z$, every term in
the numerator is divisible by $p$. If $p\nmid z$, Fermat's theorem and
the congruence of the exponents imply that all $z^{d_i}$ are congruent
modulo $p$, so their sum is again divisible by $p$.

We have $P(0)=0$ and $P(1)=1$. Also $P(2)$ is even. If $p=2$, every
term in the numerator is divisible by $4$, while if $p$ is odd, the
numerator is even and division by the odd integer $p$ preserves evenness.
Thus $h=2$. Each monomial $z^{d_i}$ has strictly increasing consecutive
gaps, and hence so does $P(z)$.

Apply \cref{cor:gap} with $a=pc$, where $c\ge1$. Since
$pc\ge2=h$ and the gaps strictly increase, the gap condition holds. Therefore
\[
\boxed{
\mathcal{R}_2\left(
x-y=c\sum_{i=1}^{p}z^{d_i}
\right)
=
c\sum_{i=1}^{p}
\left((pc+1)^{d_i}+2^{d_i}\right)+1.
}
\]

For $p=2$, the congruence condition on the exponents is vacuous. Thus,
for $m,n\ge2$,
\[
\mathcal{R}_2(x-y=c(z^m+z^n))
=
c(2c+1)^m+c(2c+1)^n+c2^m+c2^n+1.
\]
For $p=3$, the exponents need only have the same parity. Hence, whenever
$r,s,t\ge2$ are all even or all odd,
\[
\begin{aligned}
\mathcal{R}_2(x-y=c(z^r+z^s+z^t))
={}&c(3c+1)^r+c(3c+1)^s+c(3c+1)^t\\
&+c2^r+c2^s+c2^t+1.
\end{aligned}
\]
Thus the same construction gives exact Rado numbers for sums of a prime
number of powers satisfying the congruence condition above.
\subsection{Sums of consecutive powers}

Let
\[
P(z)=\sum_{r=1}^{z}r^d,
\qquad d\ge1.
\]
By Faulhaber's formula, this is an integer-valued polynomial in $z$.
We have $P(0)=0$ and $P(1)=1$, while $P(2)=1+2^d$ is odd and
$P(3)=1+2^d+3^d$ is even. Thus $h=3$ and
$P(h)=1+2^d+3^d$.

The consecutive gaps are
\[
P(n+1)-P(n)=(n+1)^d,
\]
and strictly increase. For $a\ge3$,
\[
P(a+1)-P(a)=(a+1)^d
\ge4^d
\ge2^d+3^d-1
=P(h)-2.
\]
Hence \cref{cor:gap} gives
\[
\boxed{
\mathcal{R}_2\left(
x-y=a\sum_{r=1}^{z}r^d
\right)
=
a\sum_{r=1}^{a+1}r^d
+a(1+2^d+3^d)+1
\qquad(a\ge3).
}
\]

\subsection{Polygonal numbers}

The theorem also gives a single family containing the triangular,
square, pentagonal, hexagonal, and higher polygonal numbers. For
$s\ge3$, let
\[
P(z)=\frac{(s-2)z^2-(s-4)z}{2}.
\]
The numerator is always even, so $P(z)$ is integer-valued. Also
$P(0)=0$ and $P(1)=1$. Its consecutive gaps are
\[
P(n+1)-P(n)=(s-2)n+1,
\]
and hence strictly increase. Now $P(2)=s$. If $s$ is even, then $h=2$ and $P(h)=s$. If $s$ is odd,
then $P(2)$ is odd while $P(3)=3(s-1)$ is even, so $h=3$ and
$P(h)=3(s-1)$.

If $s$ is even and $a\ge2$, then
\[
P(a+1)-P(a)=(s-2)a+1\ge s-2=P(h)-2.
\]
If $s$ is odd and $a\ge3$, then
\[
P(a+1)-P(a)
\ge3(s-2)+1
=3(s-1)-2
=P(h)-2.
\]
Thus \cref{cor:gap} gives
\[
\boxed{
\begin{aligned}
\mathcal{R}_2\left(
x-y=a\frac{(s-2)z^2-(s-4)z}{2}
\right)
&=
\frac{a}{2}\left((s-2)(a+1)^2-(s-4)(a+1)\right)\\
&\qquad+
\begin{cases}
as+1, & s\text{ even},\\
3a(s-1)+1, & s\text{ odd}.
\end{cases}
\end{aligned}
}
\]
This holds for $a\ge2$ when $s$ is even and $a\ge3$ when $s$ is odd. This one formula contains exact Rado numbers associated with every
classical polygonal number sequence.

\subsection{Simplex numbers}

For $k\ge2$, let
\[
P(z)=\binom{z+k-1}{k}.
\]
Put $q=2^{\nu_2(k)}$, where $\nu_2(k)$ is the largest nonnegative integer $r$ such that
$2^r$ divides $k$. By Kummer's theorem in the case $p=2$,
$\binom{k+t}{k}$ is odd exactly when adding $k$ and $t$ in binary
creates no carry \cite{Granville}. Since $q=2^{\nu_2(k)}$ corresponds to the least nonzero binary digit of
$k$, no carry occurs for $1\le t<q$, while adding $t=q$ produces a carry.
Hence
\[
h=1+q,
\qquad
P(h)=\binom{k+q}{k}.
\]
The consecutive gaps are
\[
P(n+1)-P(n)=\binom{n+k-1}{k-1},
\]
and strictly increase. Therefore \cref{cor:gap} applies whenever
\[
a\ge1+q
\qquad\text{and}\qquad
\binom{a+k-1}{k-1}
\ge
\binom{k+q}{k}-2.
\]
Under these conditions,
\[
\boxed{
\mathcal{R}_2\left(
x-y=a\binom{z+k-1}{k}
\right)
=
a\binom{a+k}{k}
+a\binom{k+q}{k}+1.
}
\]
The following example shows that the gap condition in \cref{cor:gap} is not necessary. Take $k=4$ and $a=5$. Then
\[
h=5,\qquad
P(5)=70,\qquad
P(6)=126,\qquad
P(7)=210,
\]
so
\[
P(6)-P(5)=56<68=P(h)-2,
\qquad
P(7)-P(6)=84.
\]
\cref{cor:gap} does not apply. However, the consecutive gaps are
strictly increasing, so \cref{cor:convex} gives the range condition.
Moreover, the nearest values to $P(6)=126$ are at distances $56$ and
$84$, so no $P(j)$ differs from $126$ by an odd integer in
$\{1,3,\ldots,67\}$. Hence the exact condition in \cref{thm:main}
holds, and
\[
\mathcal{R}_2\left(x-y=5\binom{z+3}{4}\right)=981.
\]
\subsection{Rising factorials}

The identity
\[
k!\binom{z+k-1}{k}
=
z(z+1)\cdots(z+k-1)
\]
allows the simplex number result to be written as an exact formula for
rising factorials. Apply the simplex number formula with $a=ck!$.

For $k\ge3$ and $c\ge1$, put $a=ck!$. Since $q=2^{\nu_2(k)}\le k$, we have $a\ge1+q=h$. For $k=3$, we have $q=1$, and the gap condition is immediate since
$\binom82=28>\binom43-2=2$. For $k\ge4$,
\[
\binom{ck!+k-1}{k-1}
\ge\binom{k!+k-1}{k-1}
>
\left(\frac{k!}{k-1}\right)^{k-1}
\ge8^{k-1}>4^k>\binom{2k}{k}
\ge\binom{k+q}{k}.
\]
Thus the gap condition holds. When $k=2$ and $c\ge2$, we have $a=2c\ge4>h=3$ and $P(a+1)-P(a)=a+1\ge5>P(h)-2=4$. Hence
\[
\boxed{
\begin{aligned}
&\mathcal{R}_2\left(
x-y=c\,z(z+1)\cdots(z+k-1)
\right)\\
&\qquad=
c\prod_{j=1}^{k}(ck!+j)
+
c\prod_{j=1}^{k}(q+j)
+1,
\qquad q=2^{\nu_2(k)}.
\end{aligned}
}
\]
This holds for $k\ge3$, $c\ge1$, and also for $k=2$, $c\ge2$.
The case $k=2,c=1$ lies outside the hypotheses used here.

\section{Relating Polynomial Van der Waerden and Rado numbers}

This paper is motivated by the question of when Polynomial Van der Waerden numbers and nonlinear Rado numbers, two seemingly different Ramsey quantities, are linked by an exact formula.
The connection with Polynomial Van der Waerden numbers is easiest to
see by comparing the two definitions. The number $W(aP(z),2)$ forces
two equally colored integers $x>y$ such that
\[
x-y=aP(t)
\]
for some positive integer $t$. For $(x,y,t)$ to be a monochromatic
solution of
\[
x-y=aP(z),
\]
we only need one additional fact: $t$ must have the same color as
$x$ and $y$.

For the polynomials considered here,
\[
W(aP(z),2)=aP(h)+1.
\]
The proof of \cref{thm:main} shows how this extra color condition on
$t$ is forced. In the upper bound, the same distances $aP(1)=a$ and
$aP(h)$ that determine $W(aP(z),2)$ are used to force the colors of
the possible roots. In the lower bound, the alternating coloring that
avoids these polynomial differences is used at the two ends of the
constructed coloring.

The remaining term $aP(a+1)$ records how far apart these two ends are
placed. The condition in \cref{thm:main} is exactly what prevents this
placement from creating an additional monochromatic polynomial
difference. This gives
\[
\mathcal{R}_2(x-y=aP(z))
=
aP(a+1)+W(aP(z),2).
\]

\end{document}